\documentclass[sts]{imsart}
\usepackage[centertags]{amsmath}
\RequirePackage{amsthm,amsfonts,amssymb}
\RequirePackage[authoryear]{natbib}
\RequirePackage[colorlinks,citecolor=blue,urlcolor=blue]{hyperref}
\RequirePackage{graphicx}
\usepackage{constants, bbm}
\usepackage{enumitem}
\usepackage{prodint}

\usepackage[mathscr]{euscript}
\usepackage{pifont}
\usepackage{mathtools}
\mathtoolsset{showonlyrefs}
\usepackage[T1]{fontenc}
\usepackage{fix-cm}

\usepackage{float}

\usepackage{color}
\definecolor{darkred}{RGB}{100,0,0}
\definecolor{darkgreen}{RGB}{0,100,0}
\definecolor{darkblue}{RGB}{0,0,150}
\usepackage{mathrsfs}
\usepackage{hyperref}
\hypersetup{colorlinks=true, linkcolor=darkred, citecolor=darkgreen, urlcolor=darkblue}

\startlocaldefs
\theoremstyle{plain}

\newtheorem{theorem}{Theorem}
\newtheorem{lem}{Lemma}
\newtheorem{cor}{Corollary}
\newtheorem{prp}{Proposition}

\theoremstyle{remark}
\newtheorem{Def}{Definition}
\newtheorem{example}{Example}
\newtheorem*{example*}{Example}

\def\sF{\mathscr{F}}

\def\beq{\begin{equation}} 
\def\eeq{\end{equation}}
\def\beqn{\begin{eqnarray*}}
\def\eeqn{\end{eqnarray*}}
\def\Bitem{\begin{itemize}\setlength{\itemsep}{.2in}}
\def\bitem{\begin{itemize}\setlength{\itemsep}{.05in}}
\def\eitem{\end{itemize}}
\def\Benum{\begin{enumerate}\setlength{\itemsep}{.2in}}
\def\benum{\begin{enumerate}\setlength{\itemsep}{.05in}}
\def\eenum{\end{enumerate}}
\def\bmult{\begin{multline*}}
\def\emult{\end{multline*}}
\def\bcenter{\begin{center}}
\def\ecenter{\end{center}}
\def\bframe{\begin{frame}}
\def\eframe{\end{frame}}

\def\cE{\mathcal{E}}

\def\cM{\mathcal{M}}
\def\cN{\mathcal{N}}

\def\bbR{\mathbb{R}}

\newcommand{\E}{\operatorname{\mathbb{E}}}

\def\eps{\varepsilon}

\def\1{\mathbbm{1}}

\def\sF{\mathscr{F}}

\endlocaldefs

\begin{document}

\begin{frontmatter}
\title{Deepening the Understanding of Double Robustness Geometrically}
\runtitle{Deepening the Understanding of Double Robustness Geometrically}

\begin{aug}
\author{\fnms{Andrew} \snm{Ying}\ead[label=e1]{}}


\end{aug}

\begin{abstract}
Double robustness (DR) is a widely-used property of estimators that provides protection against model misspecification and slow convergence of nuisance functions. Despite its widespread application, the theoretical foundation of DR remains underexplored. While DR is a property of global invariance along both nuisance directions, it is often implied by influence curves (ICs), which only have zero first-order derivatives in those directions locally. On the other hand, some literature proved the absence of DR estimating functions for the same estimand, under one parameterization yet was able to find one under another parameterization, highlighting the nuances in parameterization.

In this short communication, we address two key questions:
(1) Why do ICs frequently imply DR ``for free''? (2) Under what conditions would a given statistical model and parameterization support or prevent the existence of DR estimators? Using tools from semiparametric theory, we show that convexity is the crucial property that enables influence curves to imply DR. We then derive necessary and sufficient conditions for the existence of DR estimators.

Our main contribution also lies in the novel geometric interpretation of DR using information geometry, a discipline devoted to integrating global differential geometry with statistical analysis. By leveraging concepts such as parallel transport, m-flatness, and m-curvature freeness, we characterize DR in terms of invariance along submanifolds. This geometric perspective deepens the understanding of when and why DR estimators exist.

The results not only resolve apparent mysteries surrounding DR but also have practical implications for the construction and analysis of DR estimators. The geometric insights open up new connections and directions for future research. Our findings aim to solidify the theoretical foundations of a fundamental concept and contribute to the broader understanding of robust estimation in statistics.
\end{abstract}



\end{frontmatter}



\section{Introduction}

Double robustness (DR) has emerged as a critical property of estimators in various fields including causal inference \citep{robins1997causal, robins1997marginal, robins1999association, van2003unified, bang2005doubly, tchetgen2010doubly, chernozhukov2018double, ghassami2022minimax, ying2021proximal}, missing data \citep{robins1995analysis, tsiatis2006semiparametric, kang2007demystifying}, and survival analysis \citep{robins2005inverse, ying2022proximal, wang2022doubly}. In its original form, an estimator is called ``doubly robust'' if it remains consistent provided that at least one of two nuisance functions it relies is consistently estimated. This property offers protection against both model misspecification and slow convergence of nuisance functions. Existing research has primarily focused on constructing DR estimators for various estimands, employing methods such as targeted learning \citep{van2006targeted}, one-step estimator \citep{rotnitzky2021characterization}, or estimating equation estimator \citep{ghassami2022minimax}.  

Despite its widespread application, the theoretical foundation of DR remains underexplored. For example, various works \citep{chernozhukov2018double, hines2022demystifying} have used the so-called influence curves (ICs) to derive DR estimators. ICs possess the local property of having zero first-order derivatives along nuisance directions, whereas DR enjoys the stronger property of global invariance along both nuisance directions. On the other hand, \citet{tchetgen2010doubly} proved the absence of DR estimators for the same estimand, under one parameterization yet was able to find one under another parameterization, highlighting the nuances in parameterization. Therefore, in this paper, we aim to answer: 
\begin{enumerate}
    \item Why do ICs often imply DR ``for free''?
    \item Under what conditions would a given statistical model and parameterization support or prevent the existence of DR estimators?
\end{enumerate}

Addressing these questions is essential for advancing our understanding of DR and its practical application. Existing foundational work on DR \citep {robins2001comment, robins2008higher, rotnitzky2021characterization} has not fully resolved the apparent discrepancy between local and global properties, nor has it provided a complete characterization of the conditions for DR estimator existence.

In this paper, we address these questions by employing semiparametric theory, which applies geometric tools to analyze the local properties of spaces of probability distributions. Our contributions are two-fold. First, we show that convexity is the key property enabling influence curves to imply DR, thus resolving the discrepancy between local and global properties. Second, building on this insight, we provide necessary and sufficient conditions for the existence of DR estimators. To reinforce the geometric perspective, we draw on concepts from information geometry—such as parallel transport, m-flatness, and m-curvature freeness—that explore global properties in spaces of probability distributions. These tools allow us to characterize DR in terms of invariance along submanifolds, providing new insights into the existence and nature of DR estimators. This geometric interpretation not only deepens our understanding of DR but also forges connections with other fields and suggests new research directions. Our findings carry both theoretical and practical implications. By clarifying conditions for existence, our work can guide the construction of DR estimators and the selection of parameterizations. The geometric perspective offers a fresh view of DR and its relationship with other statistical concepts, solidifying the theoretical foundation of robust estimation and its applications across fields.

The remainder of the paper is organized as follows. Section \ref{sec:prep} introduces the necessary background on DR and semiparametric theory, setting the stage for our main results. Sections \ref{sec:first} and \ref{sec:second} present our findings on the role of convexity and the existence of DR estimators, respectively. Section \ref{sec:geom} develops the geometric interpretation using information geometry, providing new insights into the nature of DR. Finally, Section \ref{sec:dis} discusses the implications of our work and potential directions for future research.

\section{Background and Problem Setup}\label{sec:prep}
Define a sample space $\Omega$, $\sigma$-algebra $\sF$. We define a random vector $X$ and an associated probability density $p$ with respect to some reference measure on $\Omega$. We define $\cM = \{p\}$ as a model, that is, the collection of allowed probability densities of $X$. We use $p$ and $\E$ as a probability law and its expectation. We also define $L_0^2(p)$ as the space of mean zero square integrable functions of $X$, that is, $\{f(X): \E[f(X)] = 0, \E[f(X)^2] < \infty\}$, which is equipped with the inner product defined as $\langle f, g\rangle_p:= \E[f(X)g(X)]$. We are interested in inferring a finite-dimensional parameter $\theta = \theta(p): \cM \to \Theta \subset \bbR^k$. As shown in the examples in Section \ref{sec:exm}, our postulation of the estimand and nuisance function can accommodate both classical statistical problems with semiparametric models and recent literature trends on model-free and assumption-lean statistics \citep{chernozhukov2018double, berk2021assumption, hines2022demystifying, vansteelandt2022assumption}. Formally, given an $n$ i.i.d. sample, an estimator $\hat \theta = \hat \theta(\hat \gamma)$ might further depend on an estimate $\hat \gamma$ for a possibly infinite-dimensional nuisance function $\gamma \in \Gamma$. 

Double robustness is a desirable property of estimators that provides protection against model misspecification and slow convergence of nuisance functions when estimating scientific parameters. An estimator that relies on two nuisance functions $\gamma = (\gamma_1, \gamma_2) \in \Gamma_1 \times \Gamma_2$ is called doubly robust if $\hat \theta(\hat \gamma_1, \hat \gamma_2)$ is consistent for $\theta$ provided that either $\hat \gamma_1 - \gamma_1$ or $\hat \gamma_2 - \gamma_2$ converges to zero in some metric, where $\hat \gamma_1$ and $\hat \gamma_2$ are some sample-based functions. 

For an smooth estimand $\theta$ that has a regular and asymptotically linear estimator, DR can be further categorized into ``model double robustness'' and ``rate double robustness'' \citep{chernozhukov2018double, smucler2019unifying, rotnitzky2021characterization, ying2023asymptotic}. Model double robustness refers to estimators that are asymptotically normal when either of the parametric models for the nuisance functions $(\hat \gamma_1, \hat \gamma_2)$ is correctly specified. Rate double robustness, on the other hand, refers to estimators that are asymptotically normal when the product of the error rates of the nuisance estimators $(\hat \gamma_1, \hat \gamma_2)$ converges faster than root-$n$, even if they are estimated nonparametrically.

Every regular and asymptotically linear estimator is equivalent to an estimating function, up to some regularity conditions. See \citet{bickel1993efficient, van2000asymptotic, tsiatis2006semiparametric} for more details. To sharp the focus, we hence concentrate on population-level double robustness below. We formally define adaptive estimating function and double robustness as follows:
\begin{Def}[Adaptive estimating function]\label{def:ee}
A possibly vector function $\phi(X; \theta, \gamma)$ of $X$, parameter of interest $\theta(p)$, and possibly infinite-dimensional function $\gamma(p)$ is called an adaptive estimating function of $\theta$ when it satisfies, for any $p \in \cM$,
\begin{equation}
    \E[\phi(X; \theta(p), \gamma(p))] = 0,
\end{equation}
\begin{equation}
    \E[\phi(X; \theta', \gamma(p))] \neq 0, 
\end{equation}
for all $\theta' \neq \theta(p)$ in a neighborhood of $\theta(p)$, and
\begin{equation}
    \E[\phi(X; \theta(p), \gamma(p))^2] < \infty.
\end{equation}
We write $\phi_p(X) = \phi(X; \theta(p), \gamma(p))$ for simplicity.
\end{Def}
We call $\gamma$ the ``nuisance function'' of the estimating function. Here $\gamma(p)$ is a function from $\cM$ onto the nuisance function space. We call $(\theta(), \gamma())$ a parameterization over $\cM$. Intuitively a parameterization acts as longitude and latitude over the statistical manifold. Remark that traditionally an (non-adaptive) estimating function $\phi(X; \theta)$ does not depend on any nuisance function $\gamma$. 

Since in this paper, we focus on DR, all nuisance functions considered are in the form of $\gamma = (\gamma_1, \gamma_2)$, unless noted otherwise. Results on multiple robustness can be similarly achieved without much additional effort. We now define a population-level DR over an adaptive estimating function:
\begin{Def}[Doubly robust estimating function]\label{def:dr}
An adaptive estimating function $\phi(X; \theta, \gamma_1, \gamma_2)$ is called doubly robust if it remains unbiased provided that either one of the nuisance functions is fixed at the truth, that is,
\begin{equation}
    \E[\phi(X; \theta(p), \gamma_1(p), \gamma_2)] = \E[\phi(X; \theta(p), \gamma_1, \gamma_2(p))] = 0,
\end{equation}
for any $\gamma_1$ and $\gamma_2$.
\end{Def}
We maybe drop ``adaptive'' here because ``doubly robust'' implies that it relies on nuisance functions. Note that DR is global property of an adaptive estimating function as it allows one of $(\gamma_1, \gamma_2)$ to vary freely provided the other one is fixed at the truth. This further implies a global property on the parameterization as well along $\cM$ and how this differs from its counterparts for estimators. A DR estimating function is usually the foundation of building a DR estimator \citep{hines2022demystifying, ying2023asymptotic}, which however depends further on how one estimates the nuisance function $\gamma$. 

\subsection{Examples}\label{sec:exm}
The following examples illustrate the concept of DR estimating functions in various settings:

\begin{example}[Partially linear model]\label{exm:plm}
Consider the following partially linear model \citep{hardle2000partially}
\begin{equation}
    Y = \theta^\top A + \omega(L) + \eps,
\end{equation}
where $X = (Y, A, L)$, $\omega$ is unknown and $\E(\eps|A, L) = 0$. $A$ and $L$ are both explanatory variables whilst $A$ has linear effect and $L$ is nonlinear. The model is semiparametric since it contains both parametric component $\theta$ and nonparametric component $\omega()$. With appropriate causal conditions \citep{robins2001comment, vansteelandt2014structural}, this model is also known as the structural mean model and $\theta$ can also be understood causally. The estimating function
\begin{equation}
    \phi(X; \theta, \gamma_1, \gamma_2) = [d(A, L) - \gamma_1(X)][Y - \theta^\top A - \gamma_2(X)]
\end{equation}
is doubly robust, where
\begin{equation}
    \gamma_1(p)(X) = \E[d(A, L)|L],
\end{equation}
$d(A, L)$ is some known function of $(A, L)$, and
\begin{equation}
    \gamma_2(p)(X) = \omega(L).
\end{equation}
It is easy to show that $\phi(X; \theta, \gamma_1, \gamma_2)$ is doubly robust.
\end{example}

We then proceed to the semiparametric odds ratio model to show how DR can be absent under one parameterization but present under a different parameterization.
\begin{example}[Odds ratio model, canonical parameterization]\label{exm:odds1}
Suppose we observe $X = (Y, A, L)$, a binary $Y$, binary $A$, and some covariates $L$. Consider the semiparametric conditional odds ratio model:
\begin{equation}
    \frac{p(Y|A, L)p(y_0|a_0, L)}{p(Y|a_0, L)p(y_0|A, L)} = \psi(Y, A, L; \theta).
\end{equation}
for some baseline point $y_0$, $a_0$, where $\psi$ is a known function. We are interested in inferring $\theta$. Suppose the parameterization is given by the canonical density decomposition
\begin{equation}
    \gamma_1(p)(X) = p(A|L),
\end{equation}
and
\begin{equation}
    \gamma_2(p)(X) = p(Y|A, L).
\end{equation}
\citet{robins2001comment} has shown that there does not exist any DR estimating functions under this parameterization.
\end{example}

\begin{example}[Odds ratio model, another parameterization]\label{exm:odds2}
Continuing Example \ref{exm:odds1}, however, consider a different parameterization by
\begin{equation}
    \gamma_1(p)(X) = p(Y|a_0, L),
\end{equation}
and
\begin{equation}
    \gamma_2(p)(X) = p(A|y_0, L).
\end{equation}
\cite{chen2007semiparametric} has shown that all influence curves (which however do not have explicit form) are DR in this case. See more discussion in \citet{tchetgen2010doubly} as well.
\end{example}

\begin{example}[Average treatment effect]\label{exm:ate}
Suppose we observe $X = (Y, A, L)$ where $Y$ is an outcome of interest, $A$ is the treatment, and $L$ are baseline covariates that ensure there is no unmeasured confounding. We are interested in the average treatment mean on one treatment arm
\begin{equation}
    \theta = \E[\E(Y|A = a, L)].
\end{equation}
The estimating function
\begin{align}
    \phi(X; \theta, \gamma_1, \gamma_2)&=\gamma_2(X)(Y - \theta)\\
    &-\gamma_2(X)[\gamma_1(X) - \theta] \\
    &+ \int\gamma_1(X)d\mathbbm{1}(A = a) - \theta,
\end{align}
is doubly robust, where
\begin{equation}
    \gamma_1(p)(X) = \E(Y|A, L), 
\end{equation}
and
\begin{equation}
    \gamma_2(p)(X) = \frac{\mathbbm{1}(A = a)}{p(A|L)}.
\end{equation}
Hence
\begin{equation}
\int\gamma_1(X)d\mathbbm{1}(A = a) = \E(Y|A = a, L).
\end{equation}
See \citet{hernan2020causal} for details.
\end{example}

For readers interested in more examples, we refer to \citet{chernozhukov2018double, rotnitzky2021characterization, hines2022demystifying, ghassami2022minimax, ying2023asymptotic}.

\subsection{Semiparametric theory and influence curves}

Semiparametric theory \citep{newey1990semiparametric, bickel1993efficient, van2000asymptotic, bickel2001inference, van2003unified, tsiatis2006semiparametric, kosorok2008introduction, kennedy2017semiparametric} focuses on local behavior of estimands, the probability measure space, and estimators geometrically. This theory emerges from borrowing ``local'' concepts of differential geometry to the probability measure ``manifold.'' 

In the context of a differentiable manifold with a differentiable real-valued function defined on it, the differential is understood as a linear approximation to the function at points along the tangent vectors. These tangent vectors represent an equivalence class of differentiable curves, defined through the equivalence relation of first-order contact between the curves. To calculate the differential, one can use any direction, or equivalently, any tangent vector, to determine the directional derivative. 

In the case of a probability measure ``manifold,'' these concepts take on explicit representations. For any smooth curve on $\cM$, also known as a (one-dimensional) parametric submodel $\{p_t\}$, that passes through $p$ ($p_0 = p$), one can derive a corresponding tangent vector $s(x) = \frac{d}{dt}\big|_{t = 0}\log(p_t(x))$, which is assumed to be in $L_0^2(p)$. The tangent space is defined as the completion of the linear spans of all tangent vectors, which forms a subspace of $L_0^2(p)$, defined as below.
\begin{Def}[Tangent space]\label{def:ts}
The tangent space for an arbitrary model $\cM$ at $p$ is defined as
\begin{align*}
    &T_{p}(\cM) \\
    &=\overline{\operatorname{span}\Bigg(\Bigg\{s(X) = \frac{d}{dt}\bigg|_{t = 0}\log(p_t(X)) \in L_0^2(p):}\\
    &~~~~~~\overline{p_t \in \cM, p_0 = p\Bigg\}\Bigg)},
\end{align*}
where $\operatorname{span}(\cdot)$ is the linear span of a set and overline represents the closure of a set in $L_0^2(p)$. 
\end{Def}

\begin{figure}[t!]
\centering
\includegraphics[scale = 0.24]{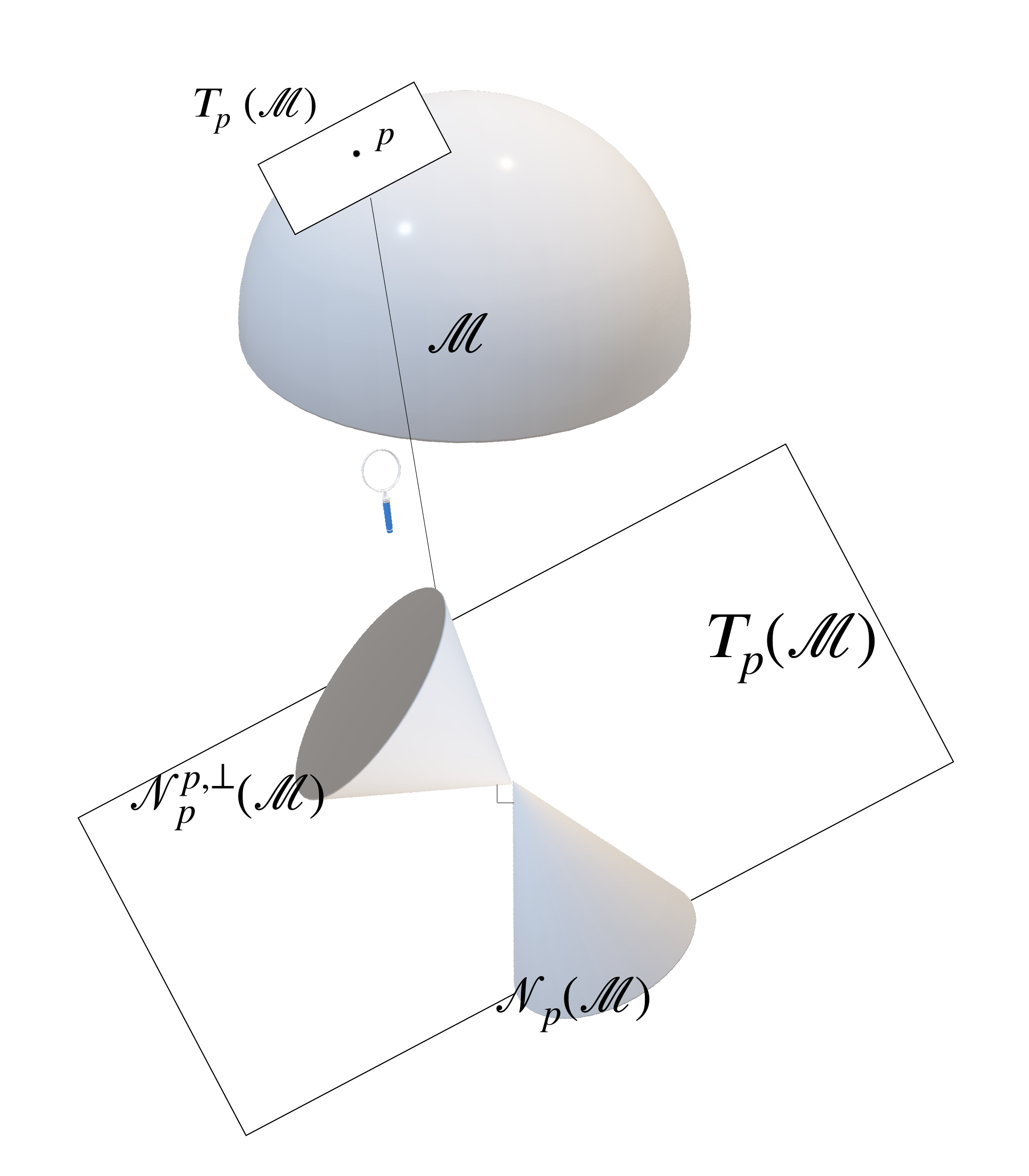}
\caption{The tangent space $T_{p}(\cM)$, the nuisance tangent space $\cN_{p}(\cM)$, and the orthogonal complement $\cN_{p}^{p, \perp}(\cM)$. Note that $\cN_{p}(\cM) \subset T_{p}(\cM)$.}
\label{fig:tangentspace}
\end{figure}
Understanding the tangent space enables one to compute the differential of the parameter of interest, which can, by Riesz's representation theorem, be represented by an influence curve. Following \citet{van2000asymptotic}, we define:
\begin{Def}[Influence curve]\label{def:if}
When $\theta: \cM \to \bbR^k$ is differentiable in the sense that there exists a continuous linear map (the differential map): $\dot \theta_p: L_0^2(p) \to \bbR^k$ such that for any parametric submodel $\{p_t\}$ with $p_0 = p$ and its corresponding score $s(X)$
\begin{equation}
    \frac{d\theta(p_t)}{dt}\bigg|_{t = 0} = \dot \theta_p(s).
\end{equation}
An influence curve $\text{IC}(p)$, being an $\bbR^k$-valued function with each component belonging to $L_0^2(p)$, is a Riesz representer of the functional $\dot \theta_p$, that is, 
\begin{equation}
    \dot \theta_p(s)= \E[\text{IC}(p)s(X)].
\end{equation}
\end{Def}
In this paper, influence curves and adaptive estimating functions (including double robust estimating functions) are all $\bbR^k$-valued functions. Any results that these functions belong to some spaces hold component-wise because all function spaces are defined for real-valued functions. 

The term ``influence curve'' is used in \citet{van2003unified, kennedy2017semiparametric}. In particular, \citet{kennedy2017semiparametric} uses it as the geometric object defined in Definition \ref{def:if}, to distinguish with the commonly used ``influence function'' of a regular and asymptotically linear estimator, though they usually coincide with each other under minor conditions. Understanding the ICs of the estimand is crucial in practice, since it often suggests viable estimators in numerous scenarios that can substantially mitigate first-order bias \citep{hines2022demystifying}. The efficient influence curve (EIC) is the unique Rieze representer when the domain of the continuous linear functional $\dot \theta_p$ is restricted to the tangent space.

Another important concept is the nuisance tangent space, a subset of the tangent space $T_p(\cM)$ where $\theta(p)$ ``does not change'' in the first order.
\begin{Def}[Nuisance tangent space]\label{def:nts}
The nuisance tangent space of $\theta(p)$ at $p$ is defined
\begin{equation}
\cN_{p}(\cM) = \overline{\operatorname{span}\left(\left\{s(X) \in T_{p}(\cM): \frac{\text{d} \theta(p_t)}{\text{d} t}\bigg|_{t = 0} = 0\right\}\right)}.
\end{equation}
\end{Def}
This way of definition is adopted in \citet[Definition 1.1]{van2003unified}, without the needs of defining a one-to-one parameterization over $\cM$ like in \citet{van2000asymptotic, tsiatis2006semiparametric}, and hence is more tailored for this draft. See Figure \ref{fig:tangentspace} for a geometric view of the tangent space $T_{p}(\cM)$, the nuisance tangent space $\cN_{p}(\cM)$, and its orthogonal complement $\cN_{p}^{p, \perp}(\cM) \subset L_0^2(\cM)$ at $p$ for the model $\cM$. In fact, it is well known that IC not only lies in but also spans the orthogonal complement of the nuisance tangent space, that is, $\overline{\operatorname{span}(\{\text{all ICs}\})} = \cN_{p}^{p, \perp}(\cM)$ \citep{van2003unified}. In fact, Definition \ref{def:if} \citep{hines2022demystifying, ying2021proximal} and computing $\cN_{p}^{p, \perp}(\cM)$ \citep{tsiatis2006semiparametric} provide two common ways of finding ICs.

We informally outline what we meant in the introduction that ICs have zero first order derivative to nuisance directions $\gamma()$. For any parametric submodel $\{p_t\}$ with the corresponding score $s(X) \in \cN_p(\cM)$, plugging in $p_t$ and take differentiation against $t$ yields
\begin{align}
    0&=\frac{d}{dt}\E_t[\phi_{p_t}(X)]\\
    &=\E[\phi_{p}(X)s(X)] + \frac{\partial}{\partial \theta}\E[\phi_p(X)]\frac{d\theta(p_t)}{dt}\\
    &+\frac{\partial}{\partial \gamma}\E[\phi_p(X)]\frac{d\gamma(p_t)}{dt}
\end{align}
Note that the first term in the second line is zero since $\phi_p(X)$ is an IC and $s(X) \in \cN_p(\cM)$. The second term is zero because definition of nuisance tangent space ensures that $\frac{d\theta(p_t)}{dt} = 0$. Therefore, we must have 
\begin{equation}
    \frac{\partial}{\partial \gamma}\E[\phi_p(X)] = 0,
\end{equation}
by the arbitrariness of $\frac{d\gamma(p_t)}{dt}$ from any $\{p_t\}$. Here the derivative $\frac{\partial}{\partial \gamma}$ is understood informally because $\gamma$ can be infinite-dimensional.

ICs only have the property that the gradients to $\gamma$ are zero, whilst DR implies the mean of an estimating function remains constant at zero, a much stronger result. Therefore, it is perhaps unsurprising that it will be implied by the global property of double robustness, officially stated as follows:
\begin{prp}\label{prp:drtoic}
For any doubly robust estimating function $\phi(X; \theta, \gamma_1, \gamma_2)$, $\phi_p(X)$ is orthogonal to the nuisance tangent space at any $p \in \cM$, that is, $\phi_p(X) \in \cN_p^{p, \perp}(\cM)$. 
\end{prp}
The intuition behind is that DR have ensured first order insensitivity to two nuisance directions, and in fact all, a DR estimating function relies on, hence being orthogonal to all nuisance tangent vectors. 

On the contrary, as we see in the examples, there are lots of cases when ICs, or more broadly up to a scale factor, vectors in the orthogonal complement of the nuisance tangent space are indeed DR or at least inspires DR estimating functions, see \citet{rotnitzky2021characterization, hines2022demystifying, ying2023asymptotic} for more examples and explanations.

\subsection{Problem formulation}
We now formally state the two key questions we aim to address in this paper:
\begin{enumerate}
    \item Why can ICs frequently imply DR ``for free''?
    \item Under what conditions would a given statistical model and parameterization support or prevent the existence of DR estimators?
\end{enumerate}



\section{Answering Question One: Convexity is Key}\label{sec:first}
As DR allows one to change one nuisance function freely while keeping the parameter of interest and the other nuisance function unchanged, it is intuitive and important to introduce and operate on the following concept, which captures the subset of the model where only one nuisance function can change.
\begin{Def}[Contour set]
The contour set of $\gamma_1()$ at $p$ is defined as the preimage of $\theta(), \gamma_2()$ of $\{p\}$, that is,
    \begin{equation}
    \cM_{\theta, \gamma_2}(p) = \{p^* \in \cM: \theta(p^*) = \theta(p), \gamma_2(p^*) = \gamma_2(p)\}.
\end{equation}
\end{Def}
Likewise, one can define $\cM_{\theta, \gamma_1}(p)$. Here are some quick observations:
\begin{enumerate}
    \item The truth $p$ lies in its contour set, that is, $p \in \cM_{\theta, \gamma_2}(p)$;
    \item $\gamma_1(p)$ is free to change while $\theta()$ and $\gamma_2()$ are fixed within the contour set;
    \item If $\phi(X; \theta, \gamma_1, \gamma_2)$ is doubly robust, then evaluating it at any $p^* \in \cM_{\theta, \gamma_2}(p)$ is still unbiased at $p$, that is,
    \begin{equation*}
        \E\{\phi_{p^*}(X)\} =\E[\phi(X; \theta(p), \gamma_1(p^*), \gamma_2(p))] = 0;
    \end{equation*}
    \item Contour set is a partition over the model. Any two points in the same contour set have the same contour set, that is, if $p^* \in \cM_{\theta, \gamma_2}(p)$, then $p \in \cM_{\theta, \gamma_2}(p^*)$, thus $\cM_{\theta, \gamma_2}(p) = \cM_{\theta, \gamma_2}(p^*)$.
    \item The tangent space of the contour set $\cM_{\theta, \gamma_2}(p)$ at $p$ is a subspace of the nuisance tangent space $\cN_{p}(\cM)$. In fact, for any parametric submodel $\{p_t\} \subset \cM_{\theta, \gamma_2}(p)$ with $p_0 = p$, the corresponding score $s(X)$ at $p$ lies in the nuisance tangent space $\cN_{p}(\cM)$ because $\theta(p)$ remains unchanged within $\cM_{\theta, \gamma_2}(p)$. Therefore, $T_p(\cM_{\theta, \gamma_2}(p)) \subset \cN_{p}(\cM)$.
\end{enumerate}
Likewise, all the above hold for $\cM_{\theta, \gamma_1}(p)$. 

We now present one of our main results, which unveils that convexity of contour sets is key for influence curves implying DR. We need an additional concept: we call $(\theta(), \gamma_1(), \gamma_2())$ ``variation independent'' if any combination of values $(\theta, \gamma_1, \gamma_2)$ can be achieved within the model $\cM$. This ensures that $\cM_{\theta, \gamma_j}(p)$ ($j = 1, 2$) offers enough variability of $\gamma_k(p)$ $(k \neq j)$, so that showing an estimating equation remains unbiased within $\cM_{\theta, \gamma_j}(p))$ leads to DR. In many papers, variation independence is often implicitly assumed or simply holds by construction, which is commonly left unstated. Variation independence allows each component of $(\theta(), \gamma_1(), \gamma_2())$ to vary freely without constraints imposed by the others. This makes sense because, with variation independence, we ensure that any specification errors in one nuisance do not inherently restrict or bias the other. As a result, the DR property can hold more reliably: if one of the nuisances is misspecified, the variation independence allows one to still provide an accurate estimate of the parameter of interest provided that the other nuisance is accurately estimated. This separation enhances the flexibility and reliability of DR estimators.
\begin{theorem}\label{thm:question1}
Suppose $(\theta(), \gamma_1(), \gamma_2())$ is variation independent. For any $p \in \cM$, when the contour sets $\cM_{\theta, \gamma_2}(p)$ and $\cM_{\theta, \gamma_1}(p)$ are convex at any $p \in \cM$, any influence curves are doubly robust.
\end{theorem}


The geometric intuition behind this result is that convexity of contour sets ensures that the tangent spaces of contour sets remains unchanged as one stays within contour sets, hence are the orthogonal complements of them. This invariance ensures that the ICs are invariant as well within contour sets, enabling the global DR property.

Theorem \ref{thm:question1} answers our first question, showing that convexity of contour sets, an assumption that often holds, is sufficient for influence curves to imply DR ``for free,'' In fact, it is easy to verify that in Examples \ref{exm:plm}, \ref{exm:odds2}, \ref{exm:ate}, $(\theta(), \gamma_1(), \gamma_2())$ is variation independent, and the contour sets $\cM_{\theta, \gamma_2}(p)$ and $\cM_{\theta, \gamma_1}(p)$ are convex, therefore in those examples we get double robustness ``for free'' from deriving influence curves. In Example \ref{exm:odds1}, however, $(\theta(), \gamma_1(), \gamma_2())$ is not variation independent.

\section{Answering Question Two: Examining nuisance directions globally}\label{sec:second}
Convexity of contour sets nicely explains why ICs imply DR in most cases. However, there are situations when contour sets are non-convex, see \citet{robins2001comment}. Now we investigate the necessary and sufficient conditions for the existence of DR estimators under a given parameterization.

\subsection{A necessity theorem}

We first state a key theorem that outlines a useful necessity result for DR.
\begin{theorem}\label{thm:necessary}
For any $p^* \in \cM_{\theta, \gamma_j}(p)$ $(j = 1, 2)$, any doubly robust estimating function $\phi(X; \theta, \gamma_1, \gamma_2)$ and any $s(X) \in T_{p^*}(\cM_{\theta, \gamma_j}(p^*)) = T_{p^*}(\cM_{\theta, \gamma_2}(p))$,
\begin{equation}
    \phi_p(X) \perp s(X)
\end{equation}
in $L_0^2(p^*)$. Therefore,
\begin{equation}
    \phi_p(X) \in T_{p^*}^{p^*, \perp}(\cM_{\theta, \gamma_j}(p)),
\end{equation}
where $T_{p^*}^{p^*, \perp}(\cM_{\theta, \gamma_j}(p))$ is the orthogonal complement of $T_{p^*}(\cM_{\theta, \gamma_j}(p))$ in $L_0^2(p^*)$ and it follows that
\begin{align*}
    \phi_p(X) \in
    &\left\{\bigcap_{p^* \in \cM_{\theta, \gamma_2}(p)} T_{p^*}^{p^*, \perp}(\cM_{\theta, \gamma_2}(p))\right\}\\
    &\bigcap \left\{\bigcap_{p^* \in \cM_{\theta, \gamma_1}(p)} T_{p^*}^{p^*, \perp}(\cM_{\theta, \gamma_1}(p))\right\}.
\end{align*}
\end{theorem}
This theorem can be seen as a generalization of the necessity theorem outlined by \citet[Lemma 1]{robins2001comment}. It establishes that a DR estimating function is orthogonal to the nuisance tangent space, not just at the true distribution $p$, but at any distribution $p^*$ within the contour sets of either $\gamma_1()$ or $\gamma_2()$. By showing the intersection set is empty, one might falsify the existence of a DR estimating function.  

Now we may strengthen our understanding of previous unofficial usage of ``local'' versus ``global'' more formally and mathematically. From Proposition \ref{prp:drtoic}, we can show that
\begin{equation*}
     \cN_p^{p, \perp}(\cM) = T_p^{p, \perp}(\cM_{\theta, \gamma_2}(p)) \bigcap T_p^{p, \perp}(\cM_{\theta, \gamma_1}(p)).
\end{equation*}
Note that this is local as we only compute the orthogonal complement of tangent spaces of the contour sets at $p$. However, by Theorem \ref{thm:necessary}, a DR estimating function resides in a smaller set, the intersection of $\bigcap_{p^* \in \cM_{\theta, \gamma_2}(p)} T_{p^*}^{p^*, \perp}(\cM_{\theta, \gamma_2}(p))$, together with $\bigcap_{p^* \in \cM_{\theta, \gamma_1}(p)} T_{p^*}^{p^*, \perp}(\cM_{\theta, \gamma_1}(p))$. In other words, it needs to be orthogonal to many more tangent vectors from other points in the contour sets and hence being a global concept. In fact, we show in the next subsection, with minor additional conditions on the parameterization $(\theta(), \gamma_1(), \gamma_2())$, the necessity result can be strengthened into a sufficiency theorem.

\subsection{A sufficiency theorem}

To obtain a sufficiency result in the absence of convexity, we need to introduce the following concept. We call an arbitrary model $\cM$ is ``smoothly path-connected'' if, for any pair $p$, $p^* \in \cM$, there exists a smooth curve $\{p_t\} \subset \cM$ such that $p_0 = p$ and $p_1 = p^*$. 
\begin{Def}[$\theta$-connectedness]
A parameterization $(\theta(), \gamma_1(), \gamma_2())$ is called $\theta$-connected if the contour sets $\cM_{\theta, \gamma_2}(p)$ and $\cM_{\theta, \gamma_1}(p)$ are smoothly path-connected.
\end{Def}
$\theta$-connectedness typically holds for a reasonable parameterization. Note that the convexity of both contour sets $\cM_{\theta, \gamma_2}(p)$ and $\cM_{\theta, \gamma_1}(p)$ at any $p \in \cM$ implies $\theta$-connectedness. The $\theta$-connectedness allows us to ``move'' vectors between any two distributions within a contour set, or equally, moving the orthogonal complements of tangent space of the contour sets to the same point $p$. 

We now state another key theorem that provides a sufficient condition for the existence of DR estimators.
\begin{theorem}\label{thm:iff}
Suppose that $(\theta(), \gamma_1(), \gamma_2())$ is variation independent and $\theta$-connected. An adaptive estimating function $\phi(X; \theta, \gamma_1, \gamma_2)$ is doubly robust, if for any $p \in \cM$, 
\begin{align*}
    \phi_p(X) \in 
    &\left\{\bigcap_{p^* \in \cM_{\theta, \gamma_2}(p)}T_{p^*}^{p^*, \perp}(\cM_{\theta, \gamma_2}(p))\right\} \\
    &\bigcap \left\{\bigcap_{p^* \in \cM_{\theta, \gamma_1}(p)}T_{p^*}^{p^*, \perp}(\cM_{\theta, \gamma_1}(p))\right\}.
\end{align*}
\end{theorem}

Theorems \ref{thm:necessary} and \ref{thm:iff} provide a complete characterization of the existence of DR estimators under a given parameterization. It states that if an adaptive estimating function is orthogonal to the tangent spaces of both contour sets at all distributions within the contour sets, then it is doubly robust. It generalizes the algorithm used in \citet{tsiatis2006semiparametric} by computing $\cN_{p}^{p, \perp}(\cM)$ to find ICs.

It will be interesting to apply Theorem \ref{thm:iff} to the case when contour sets are convex and see how it becomes Theorem \ref{thm:question1}. We will postpone this to Section \ref{sec:geom}, where information geometry can help to build more interpretation and intuition. A preview is that convexity ensures that orthogonal complements of tangent space of the contour sets remain invariant. Without convexity, they might change.

Theorems \ref{thm:necessary} and \ref{thm:iff} answer our second question. By checking the orthogonality condition, we can easily verify whether a given estimating function is doubly robust under a specific parameterization. Also, by checking whether the intersection is empty, one can prove the absence of a DR estimating function.

In the next section, we will explore the geometric aspects of DR in more depth, using tools from information geometry to gain further insights into the structure of the statistical model and the properties of DR estimators.

\section{Geometric understanding by information geometry}\label{sec:geom}
There is a long history of studies on geometry of manifolds of probability distributions, for instance, semiparametric theory. However, as we mentioned in the earlier section, semiparametric theory only concerns the local property of the statistical model manifold and estimand, for instance, it only operationalizes over one tangent space at one point at a time. On the other hand, DR is a global property within contour sets over the statistical manifold. Therefore, intuitively we need more advanced geometric tools to investigate relation among tangent spaces at different points. 

Information geometry \citep{amari1997information, amari2000methods, amari2016information, ay2017information} leverages more ``global'' differential geometry concepts onto the probability ``manifold'' $\cM$. Here we give some brief historical notes on how information geometry has been developed and leveraged into statistics. In 1930, Hotelling observed that a family of parametric probability density functions possesses the structure of a Riemannian manifold \citep{hotelling1930spaces}. Although his contribution was presented only as an abstract, it is widely regarded as one of the foundational works in the development of information geometry. Independently, Rao later formalized this idea by publishing a journal paper that established the Riemannian manifold structure of statistical models \citep{rao1945information}, which is a monumental work from which information geometry has emerged. Later, \citet{efron1975defining} investigated old unpublished calculations by R.A. Fisher and elucidated the results by defining the statistical curvature of a statistical model. This work was commented on by A.P. Dawid in discussions of Efron’s paper \citep{dawid1975invited}, where the e- and m-connections were suggested. In 1982, Amari studied dual affine connections, independently of the researches mentioned above, and applied them to statistics \citep{amari1982differential}. Following Efron’s, Dawid’s, and Amari's works, \citet{amari1982differential} further systematically developed the differential geometry of statistical models and elucidated its dualistic nature. It was applied to statistical inference to establish a higher-order statistical theory \citep{amari1982differential, amari1985differential, kumon1983geometrical}. A number of competent researchers have joined from the fields of statistics, vision, optimization, machine learning, etc. Many international conferences have been organized on this subject. See \citet{matsuzoe2024half} for a review of its history.

In differential geometry, affine connection and the corresponding parallel transport are tools to move tangent vector and maintain direction between tangent spaces on curved surfaces. An affine connection provides a rule for how to compare directions (vectors) at different points on a curved surface, or manifold. It tells us how a direction can ``change'' as we move across the manifold. Parallel transport is the process of moving a vector along a curve while keeping it as ``unchanged'' as possible according to the connection. This is crucial because it helps us understand how objects ``move'' on a manifold in a way that respects the curvature.

In information geometry, two types of parallel transport are defined based on two affine connections on the probability ``manifold'': the e-connection and the m-connection. The concept of e-parallel transport arises from the structure of exponential families of distributions. The e-connection is aligned with exponential parameterizations of distributions, where natural parameters linearly affect the log-probability. The m-parallel transport, on the other hand, is associated with mixture models. The m-connection corresponds to paths in which distributions are combined as convex mixtures. More importantly, they have explicit forms and their definitions do not rely on introducing the corresponding affine connection. Also, their computation does not need to go through whole curves, but only starting point and end point.

\begin{Def}[e-parallel transport]\label{def:ept}
The e-parallel transport of a vector $\phi(X)$ from $L_0^2(p)$ to $L_0^2(p^*)$ is defined as
\begin{equation}
    \prod^e_{p \to p^*}\phi(x) = \phi(x) - \E^*\{\phi(X)\},
\end{equation}
whenever it is square integrable with respect to $p^*$.
\end{Def}

\begin{figure}[t!]
    \centering
    \includegraphics[scale = 0.18]{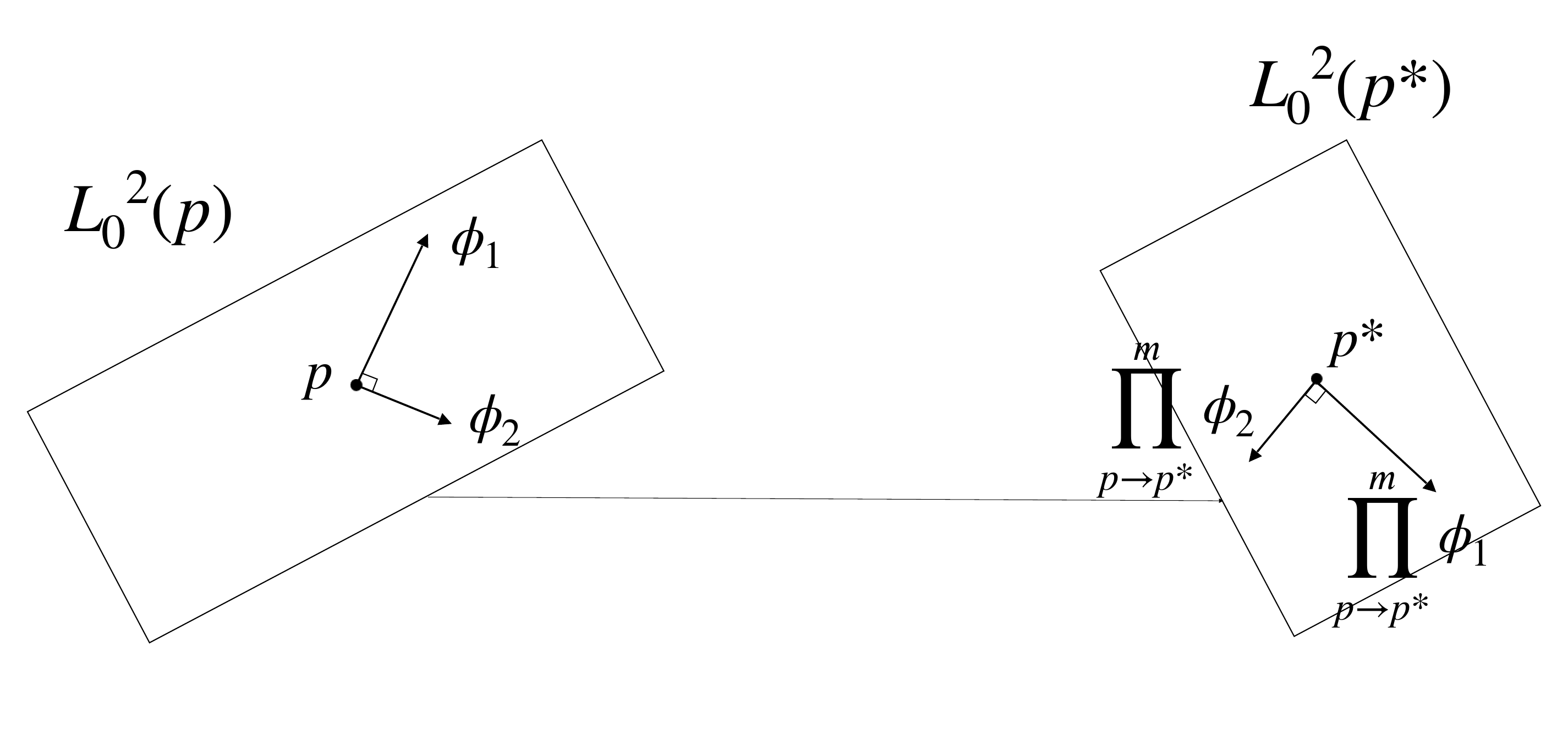}
    \caption{A graphical illustration of e- and m-parallel transports, which preserve the inner product of a pair of vectors. In particular, as shown in the graph, two perpendicular vectors ($\phi_1, \phi_2$) remain perpendicular after transportation.}
    \label{fig:emtransport}
\end{figure}

\begin{Def}[m-parallel transport]
The m-parallel transport of a vector $\phi(X)$ from $L_0^2(p)$ to $L_0^2(p^*)$ is defined as
\begin{equation}
    \prod^m_{p \to p^*}\phi(x) = \frac{p(x)}{p^*(x)}\phi(x),
\end{equation}
whenever it is square integrable with respect to $p^*$.
\end{Def}
These two parallel transports are dual with respect to the inner product, that is,
\begin{equation}
    \langle \phi_1, \phi_2\rangle_{p} = \left\langle\prod^e_{p \to p^*}\phi_1, \prod^m_{p \to p^*}\phi_2\right\rangle_{p^*}.
\end{equation}
Figure \ref{fig:emtransport} provides a visual illustration of these concepts.
We now establish a key result connecting e-parallel transport and double robustness on a definitional level.
\begin{prp}\label{prp:etransport}
Suppose $(\theta(), \gamma_1(), \gamma_2())$ is variation independent. An adaptive estimating function \\
$\phi(X; \theta, \gamma_1, \gamma_2)$ is doubly robust if and only if it remains e-parallel transport invariant along contour sets $\cM_{\theta, \gamma_2}(p)$ and $\cM_{\theta, \gamma_1}(p)$ for any $p \in \cM$.
\end{prp}

With duality between e- and m-parallel transports, and Proposition \ref{prp:etransport}, Theorems \ref{thm:necessary} and \ref{thm:iff} lead to the following corollary, characterizing DR through m-parallel tranports.
\begin{cor}[Theorem \ref{thm:necessary}]\label{cor:cor1}
An adaptive estimating function $\phi(X; \theta, \gamma_1, \gamma_2)$ is doubly robust if, for any $p \in \cM$, 
\begin{align*}
    \phi_p(X) \in 
    &\left\{\bigcap_{p^* \in \cM_{\theta, \gamma_2}(p)}\left\{\prod^m_{p^* \to p}T_{p^*}(\cM_{\theta, \gamma_2}(p))\right\}^{p, \perp}\right\} \\
    &\bigcap \left\{\bigcap_{p^* \in \cM_{\theta, \gamma_1}(p)}\left\{\prod^m_{p^* \to p}T_{p^*}(\cM_{\theta, \gamma_1}(p))\right\}^{p, \perp}\right\}.
\end{align*}  
\end{cor}

\begin{figure}[t!]
    \centering
    \includegraphics[scale = 0.18]{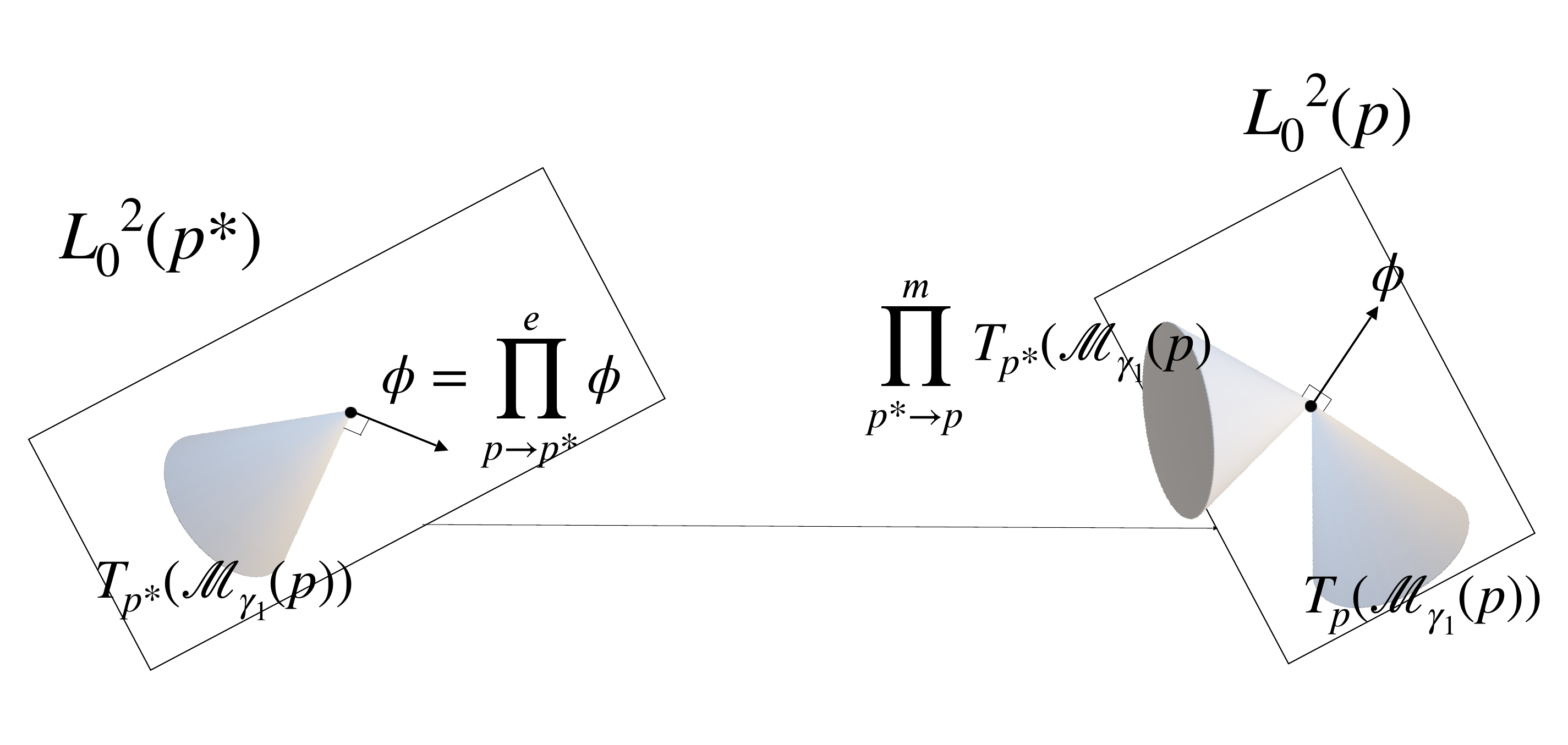}
    \caption{E-parallel transport of an vector $\phi$ in the tangent space and m-parallel transport of the tangent space of a contour set.}
    \label{fig:emtransport2}
\end{figure}

\begin{cor}[Theorem \ref{thm:iff}]\label{cor:cor2}
Suppose a parameterization $(\theta(), \gamma_1(), \gamma_2())$ is variation independent and $\theta$-connected. An adaptive estimating function $\phi(X; \theta, \gamma_1, \gamma_2)$ is doubly robust if, for any $p \in \cM$, 
\begin{align*}
    \phi_p(X) \in 
    &\left\{\bigcap_{p^* \in \cM_{\theta, \gamma_2}(p)}\left\{\prod^m_{p^* \to p}T_{p^*}(\cM_{\theta, \gamma_2}(p))\right\}^{p, \perp}\right\} \\
    &\bigcap \left\{\bigcap_{p^* \in \cM_{\theta, \gamma_1}(p)}\left\{\prod^m_{p^* \to p}T_{p^*}(\cM_{\theta, \gamma_1}(p))\right\}^{p, \perp}\right\}.
\end{align*}  
\end{cor}
See Figure \ref{fig:emtransport2} for a geometric view of the above statements.

These corollaries provide alternative necessary and sufficient conditions for DR estimators in terms of the orthogonality of the adaptive estimating function to the m-parallel transported tangent spaces of the contour sets.

\citet{amari1997information} investigated the existence of an estimating function (without nuisance functions) under a semiparametric model. They introduced the concepts of m-flatness and m-curvature freeness, which capture the behavior of the nuisance tangent spaces under m-parallel transport, to enrich geometric understanding of their topic. We remark that the theory of \citet{amari1997information} has not been mathematically justified yet, but is still useful for geometrical thinking. Following their spirits, here we generalize these concepts to understand DR.

\begin{Def}[m-flat parameterization]\label{def:mflat}
A model $\cM$ is called m-flat with respect to the parameterization $(\theta(), \gamma_1(), \gamma_2())$ if 
\begin{equation}
    \prod^m_{p^* \to p}T_{p^*}(\cM_{\theta, \gamma_j}(p)) \subset T_{p}(\cM_{\theta, \gamma_j}(p)),
\end{equation}
for any $p$ and $p^*$, and $j = 1, 2$.
\end{Def}
Intuitively, m-flatness requires that $T_{p^*}(\cM_{\theta, \gamma_j})$ remains barely changed (within $T_{p}(\cM_{\theta, \gamma_j})$) under a m-parallel transport. The convexity of $\cM_{\theta, \gamma_j}$ indeed implies m-flatness. 
\begin{lem}\label{lem:convexitytomflat}
If $\cM_{\theta, \gamma_2}(p)$ and $\cM_{\theta, \gamma_1}(p)$ are convex, then the model $\cM$ is m-flat with respect to $(\theta(), \gamma_1(), \gamma_2())$.
\end{lem}
Note that ``m-flatness'' is strictly weaker than convexity because it only requires the comparison of two tangent vectors at two points without requiring the contour sets to be connected whilst convexity implies connectedness. The weaker condition ``m-flatness'' can still achieve Theorem \ref{thm:question1}:
\begin{cor}\label{cor:cor3}
Suppose $(\theta(), \gamma_1(), \gamma_2())$ is variation independent and $\theta$-connected. If the model $\cM$ is m-flat with respect to it, any influence curves are doubly robust.
\end{cor}


This corollary is proved as follows. Note that by Theorem \ref{thm:iff} and definition of m-flatness, we have:
\begin{align*}
    \text{IC}(p) &\in\cN_p^{p, \perp}(\cM)\\
    &=T_p^{p, \perp}(\cM_{\theta, \gamma_2}(p)) \bigcap T_p^{p, \perp}(\cM_{\theta, \gamma_1}(p))\\
    &\subset \left\{\bigcap_{p^* \in \cM_{\theta, \gamma_2}(p)}\left\{\prod^m_{p^* \to p}T_{p^*}(\cM_{\theta, \gamma_2}(p))\right\}^{p, \perp}\right\} \\
    &\bigcap \left\{\bigcap_{p^* \in \cM_{\theta, \gamma_1}(p)}\left\{\prod^m_{p^* \to p}T_{p^*}(\cM_{\theta, \gamma_1}(p))\right\}^{p, \perp}\right\},
\end{align*}  
which by Corollary \ref{cor:cor2} leads to double robustness of ICs. 

In fact, Theorem \ref{thm:question1} is immediate from Lemma \ref{lem:convexitytomflat} and Corollary \ref{cor:cor3}, emphasizing the helpfulness of considering and leveraging information geometry.

Another concept from \citet{amari1997information}, which is weaker than m-flatness yet preludes double robustness, is:
\begin{Def}[m-curvature free parameterization]
A model $\cM$ is called m-curvature free with respect to the parameterization $(\theta(), \gamma_1(), \gamma_2())$ if 
\begin{equation}
    \text{EIC}(p) \in \left\{\prod^m_{p^* \to p}T_{p^*}(\cM_{\theta, \gamma_j}(p))\right\}^{p, \perp},
\end{equation}
for any $p$ and $p^*$, and $j = 1, 2$.
\end{Def}
m-curvature freeness is weaker than m-flatness as $\text{EIC}(p) \in \cN_{p}^{p, \perp}(\cM) \subset \left\{\prod^m_{p^* \to p}\cN_{p^*}(\cM_{\theta, \gamma_j}(p))\right\}^{p, \perp}$ when m-flatness holds. Intuitively m-flatness requires that $\{\prod^m_{p^* \to p}T_{p^*}(\cM_{\theta, \gamma_j}(p)\}^{p, \perp} \supset T_{p}(\cM_{\theta, \gamma_j}(p))$ whilst the m-curvature freeness weakens this argument but still requires that the most efficient direction EIC is contained in $\{\prod^m_{p^* \to p}T_{p^*}(\cM_{\theta, \gamma_j})\}^{p, \perp}$. The m-flatness ensures that the orthogonal complements of tangent space of the contour sets remain invariant within the contour sets (shown in the next subsection) while m-curvature freeness ensures that these orthogonal complements do not rotate a lot.
\begin{cor}
Suppose $(\theta(), \gamma_1(), \gamma_2())$ is variation independent and $\theta$-connected. If the model $\cM$ is m-curvature free with respect to it, the efficient influence curve is doubly robust.
\end{cor}

These corollaries highlight the role of the geometry of the tangent spaces of contour sets in determining the existence of DR estimators. In \citet{amari1997information}, nuisance tangent spaces play an important role for the geometry of estimating functions, whereas in this paper, tangent spaces of contour sets, which are subspaces of nuisance tangent spaces, play the same and essential role for the geometry of DR estimating functions. In particular, they show that certain ``nice'' geometric properties, such as m-flatness and m-curvature freeness, are sufficient for the existence of DR estimators.

The information geometric perspective developed in this section provides a deeper understanding of the global structure of the statistical model and the properties of DR estimators. By studying the behavior of adaptive estimating functions and tangent spaces of contour sets under parallel transport, we gain new insights into the geometric nature of DR.

\section{Discussion and future directions}\label{sec:dis}
In this paper, we have deepened the geometric understanding of DR by leveraging both semiparametric theory and information geometry. 

We have shown that convexity of the contour sets of the statistical model is a sufficient condition for influence curves to imply DR ``for free,'' resolving the apparent discrepancy between the local and global robustness properties of DR estimators. We have also provided necessary and sufficient conditions for the existence of DR estimators.

Furthermore, we have introduced geometric concepts, such as m-flatness and m-curvature freeness, which capture the behavior of the tangent spaces of contour sets under parallel transport. These concepts provide a deeper understanding of the global structure of the statistical model and the geometric properties that enable the existence of DR estimators.

Our findings have both theoretical and practical implications. From a theoretical perspective, our results contribute to the foundational understanding of DR and highlight the importance of geometric considerations in the study of robust estimators. The geometric characterizations we provide can serve as a basis for further investigations into the properties and behavior of DR estimators in different settings. From a practical perspective, our results can guide the construction of DR estimators and the choice of parameterizations in applied settings. The necessary and sufficient conditions we establish can be used to check the existence of DR estimators for a given statistical model and parameterization, informing the development of robust inference procedures. The geometric insights we provide can also aid in the design of more efficient and stable estimators by taking into account the structure of the statistical model.

One future direction that can be groundbreaking is: given a model $\cM$ and a parameter of interest $\theta$, if no estimating function exists but only adaptive ones, how to find a parameterization $\gamma(p)$ such that one can construct doubly (or multiply) robust estimators? One idea is to construct $\gamma = (\gamma_1, \cdots, \gamma_K)$ such that at any $p \in \cM$, in the preimage of $\theta()$, that is, $\{p^* \in \cM: \theta(p^*) = \theta(p)\}$, the further contour sets $\cM_{\theta, \gamma_j}(p)$ are all convex. Intuitively, the larger $K$, the easier this task is.



\begin{acks}[Acknowledgments]
The author would like to thank Jiangang Ying, Cong Ding for countless discussions. The author would also like to thank Eric J. Tchetgen Tchetgen for inspiration of this idea. The author is very grateful for two anonymous reviewers for greatly improving this draft. The author would like to acknowledge the helpful discussions and insights gained through interactions with the Claude AI assistant, which contributed to refining the presentation of this work but not substantive content.
\end{acks}
%

\appendix

\section{Proofs}

\subsection*{Proof of Proposition \ref{prp:drtoic}}
For any parametric submodel $\{p_t\}$ with corresponding score $s(X) \in \cN_p(\cM)$, plugging in $p_t$ and differentiating against $t$ yields
\begin{align}
    0&=\frac{d}{dt}\bigg|_{t = 0}\E_t[\phi_{p_t}(X)]\\
    &=\frac{d}{dt}\bigg|_{t = 0}\E_t[\phi(X; \theta(p_t), \gamma_1(p_t), \gamma_2(p_t))]\\
    &=\E[\phi_{p}(X)s(X)] + \frac{d}{dt}\bigg|_{t = 0}\E[\phi(X; \theta(p_t), \gamma_1(p_t), \gamma_2(p_t))]\\
    &=\E[\phi_{p}(X)s(X)] + \frac{d}{dt}\bigg|_{t = 0}\E[\phi(X; \theta(p_t), \gamma_1(p), \gamma_2(p))]\\
    &+\frac{d}{dt}\bigg|_{t = 0}\E[\phi(X; \theta(p), \gamma_1(p_t), \gamma_2(p))]\\
    &+\frac{d}{dt}\bigg|_{t = 0}\E[\phi(X; \theta(p), \gamma_1(p), \gamma_2(p_t))]
\end{align}
Note that the second term is zero because 
\begin{align}
    &\frac{d}{dt}\bigg|_{t = 0}\E[\phi(X; \theta(p_t), \gamma_1(p), \gamma_2(p))] \\
    &= \frac{d}{d\theta}\bigg|_{\theta = \theta(p)}\E[\phi(X; \theta, \gamma_1(p), \gamma_2(p))]\frac{d\theta(p_t)}{dt}\bigg|_{t = 0}
\end{align}
and the fact that the definition of nuisance tangent space implies that $\frac{d\theta(p_t)}{dt} = 0$. The third and the fourth terms are zero by double robustness. It follows that the first term is zero, implying $\phi_p(X) \perp s(X)$ at $p$ and hence the conclusion is obtained since $s(X)$ is an arbitrary score in $\cN_p(\cM)$.

\subsection*{Proof of Theorem \ref{thm:question1}}
It turns out when a model is convex, its tangent space will have an explicit form:
\begin{lem}\label{lem:convexity}
When an arbitrary model $\cM$ is convex, the tangent space for it at law $p \in \cM$ is
\begin{equation}
    T_{p}(\cM) = \overline{\operatorname{span}\left\{\frac{p^*(X)}{p(X)} - 1: p^* \in \cM\right\}}.
\end{equation}
\end{lem}
\begin{proof}
For any $p^* \in \cM$, by convexity of $\cM$, the parametric submodel $\{p_t = tp^* + (1 - t)p\}_{t \in [0, 1]} \subset \cM$, then
\begin{equation}
    \frac{p^*}{p} - 1 = \lim_{t \to 0} \frac{tp^* + (1 - t)p - p}{tp} = \frac{\textup{d}}{\textup{d}t}\bigg|_{t = 0}\log p_t
\end{equation}
is in $T_{p}(\cM)$. This has shown one direction of the lemma, that is,
\begin{equation}
    T_{p}(\cM) \supset \overline{\operatorname{span}\left\{\frac{p^*(X)}{p(X)} - 1: p^* \in \cM\right\}}.
\end{equation}
Now since $p^* \in \cM$ is arbitrary, the other direction follows.
\end{proof}

First, an immediate implication for the influence curve is
\begin{align}
    \phi_p(X) \in \cN_{p}^{p, \perp}(\cM) \subset T_{p}^{p, \perp}(\cM_{\theta, \gamma_2}(p)),
\end{align}
since $T_{p}(\cM_{\theta, \gamma_2}(p)) \subset \cN_{p}(\cM)$. That is, by the convexity of $\cM_{\theta, \gamma_2}(p)$ and Lemma \ref{lem:convexity}, for any $p^* \in \cM_{\theta, \gamma_2}(p)$,
\begin{align}
    0 &=-\E\left\{\phi_p(X)\left[\frac{p^*(X)}{p(X)} - 1\right]\right\}\\
    &=\E^*\left\{\phi_p(X)\left[\frac{p(X)}{p^*(X)} - 1\right]\right\},
\end{align}
where $\E^*$ means taking expectation with respect to $p^*$. Therefore,
\begin{equation}
    \phi_p(X) \in T_{p^*}^{p^*, \perp}(\cM_{\theta, \gamma_2}(p)).
\end{equation}
Define $g(t) = \E_t[\phi_{p^*}(X)] $, where $p_t := (1 - t)p^* + tp$ $(0 \leq t \leq 1)$ and $p \in \cM_{\theta, \gamma_2}(p^*)$. By replacing $p$ and $p^*$ with $p^*$ and $p_t \in \cM_{\theta, \gamma_2}(p^*)$ in $\phi_p(X) \in T_{p^*}^{p^*, \perp}$ respectively,
\begin{equation}
    \frac{dg(t)}{dt} = \E_t\left[\phi_{p^*}(X)s(X)\right] = 0,
\end{equation}
since $s(X)$ is the score function of $p_t$ and 
\begin{align}
    s(X) \in T_{p_t}(\cM_{\theta, \gamma_2}(p^*)). 
\end{align}
By an application of ordinary differential equation theory, we have
\begin{align}
    0 = g(0) = g(1) &= \E[\phi_{p^*}(X)]\\
    &= \E[\phi(X; \theta(p), \gamma_1(p^*), \gamma_2(p))].
\end{align}
Hence, by variation independence of $\theta(p), \gamma_1(p), \gamma_2(p)$, for any $\gamma_1^* \in \Gamma_1$, one can find $p^* \in \cM_{\theta, \gamma_2}(p)$ such that $\gamma_1(p^*) = \gamma_1^*$. Thus we have shown that
\begin{equation}
    \E[\phi(X; \theta(p), \gamma_1^*, \gamma_2(p))] = 0.
\end{equation}
In a similar way, $\E[\phi(X; \theta(p), \gamma_1(p), \gamma_2^*)] = 0$ for any $\gamma_2^* \in \Gamma_2$ can be proved. Therefore, $\phi(X; \theta, \gamma_1, \gamma_2)$ is doubly robust.

First by definition of contour sets and double robustness, for any parametric submodel $\{p_t\} \subset \cM_{\theta, \gamma_2}(p)$ with $p_0 = p$, we have
\begin{align}
    \E_t\{\phi_p(X)\}
    &=\E_t[\phi(X; \theta(p), \gamma_1(p), \gamma_2(p))]\\
    &=\E_t[\phi(X; \theta(p_t), \gamma_1(p), \gamma_2(p_t))]\\
    &=\E_t[\phi(X; \theta(p_t), \gamma_1(p_t), \gamma_2(p_t))]\\
    &=0,
\end{align}
where $\E_t$ means taking expectation with respect to $p_t$. Therefore,
\begin{align}
    0&=\frac{\textup{d}}{\textup{d}t}\bigg|_{t = 0}\E_t\{\phi_p(X)\} =\E_t[\phi(X; \theta(p), \gamma_1(p), \gamma_2(p))]\\
    &=\int\phi_p(x)\left[\frac{\textup{d}}{\textup{d}}{\textup{d}t}\bigg|_{t = 0}\log p_t(x)\right]p^*(x)\textup{d}x\\
    &=\E^*[\phi_p(X)s(X)],
\end{align}
where $s(X)$ is the score function (tangent vector) of $p_t$ at $t = 0$. Since $\{p_t\}$ is an arbitrary parametric submodel in $\cM_{\theta, \gamma_2}(p)$ with $p_0 = p^*$, this means that for any $p^* \in \cM_{\theta, \gamma_2}(p)$, any doubly robust estimating function $\phi(X; \theta, \gamma_1, \gamma_2)$ and any $s(X) \in T_{p^*}(\cM_{\theta, \gamma_2}(p^*)) = T_{p^*}(\cM_{\theta, \gamma_2}(p))$,
\begin{equation}
    \phi_p(X) \perp s(X),
\end{equation}
in $L_0^2(p^*)$. In a similar way, we can also prove this orthogonality for $s(X) \in T_{p^*}(\cM_{\theta, \gamma_1}(p^*)) = T_{p^*}(\cM_{\theta, \gamma_1}(p))$.

\subsection*{Proof of Theorem \ref{thm:iff}}
For any $p^* \in \cM_{\theta, \gamma_2}(p)$, define
\begin{equation}
    g(t) = \E_t[\phi_{p^*}(X)],
\end{equation} 
where $p_t$ $(0 \leq t \leq 1)$ is a smooth curve connecting $p = p_1$ and $p^* = p_0$. This curve can be taken in $\cM_{\theta, \gamma_2}(p) = \cM_{\theta, \gamma_2}(p^*)$ since $(\theta, \gamma_1, \gamma_2)$ is $\theta$-connected. Then by the condition for $\phi_{p^*}(X)$,
\begin{equation}
    \frac{dg(t)}{dt} = \E_t\left[\phi_{p^*}(X)s(X)\right] = 0,
\end{equation}
where $s(X)$ is the score function of $p_t$ and hence
\begin{align}
    s(X) \in T_{p_t}(\cM_{\theta, \gamma_2}(p^*)),
\end{align}
By an application of ordinary differential equation theory, we have 
\begin{align}
    0 = g(0) = g(1) &= \E[\phi_{p^*}(X)]\\
    &= \E[\phi(X; \theta(p), \gamma_1(p^*), \gamma_2(p))].
\end{align}
Again, by variation independence of $\theta(p), \gamma_1(p), \gamma_2(p)$, for any $\gamma_1^* \in \Gamma_1$, one can find $p^* \in \cM_{\theta, \gamma_2}(p)$ such that $\gamma_1(p^*) = \gamma_1^*$. Thus we have shown that
\begin{equation}
    \E[\phi(X; \theta(p), \gamma_1^*, \gamma_2(p))] = 0.
\end{equation}
In a similar way, $\E[\phi(X; \theta(p), \gamma_1(p), \gamma_2^*)] = 0$ for any $\gamma_2^* \in \Gamma_2$ can be proved.  


\subsection*{Proof of Proposition \ref{prp:etransport}}
We only need to show for $\gamma_1$. ``If'' part: for any $\gamma_1$, because of variation independence of $(\theta(), \gamma_1(), \gamma_2())$, we can find $p^* \in \cM_{\theta, \gamma_2}(p)$ such that $\gamma_1(p^*) = \gamma_1$. Also since $\phi(X; \theta, \gamma_1, \gamma_2)$ is e-parallel transport invariant along the contour set $\cM_{\theta, \gamma_2}(p)$, 
\begin{align*}
    \phi_{p^*}(X)&= \prod^e_{p^* \to p}\phi_{p^*}(X)\\
    &= \phi_{p^*}(X) - \E[\phi_{p^*}(X)],
\end{align*}
implying $\E[\phi_{p^*}(X)] = 0$. We have
\begin{align*}
    &\E[\phi(X; \theta(p), \gamma_1, \gamma_2(p))]\\
    &=\E[\phi(X; \theta(p), \gamma_1(p^*), \gamma_2(p))]\\
    &=\E[\phi_{p^*}(X)] = 0.
\end{align*}
``Only if'' part is trivial.

\subsection*{Proof of Lemma \ref{lem:convexitytomflat}}
We will prove a more general result that if a model $\cM$ is convex,
\begin{equation*}
	\prod^m_{p^* \to p}T_p^*(\cM) \subset T_p(\cM),
\end{equation*}
that is, the model $\cM$ is  m-flat in the sense defined in \citet{amari1997information}. Note from Lemma \ref{lem:convexity}, we have
\begin{equation}
    T_{p}(\cM) = \overline{\operatorname{span}\left\{\frac{p^*(X)}{p(X)} - 1: p^* \in \cM\right\}}.
\end{equation}
With the boundedness of the m-parallel transport and its linearity, it suffices to show that
\begin{align}
    &\prod^m_{p^* \to p}\left\{\frac{p^{**}(X)}{p^*(X)} - 1: p^{**} \in \cM\right\}\\
    &=\left\{\frac{p^{**}(X) - p^{*}(X)}{p(X)}: p^{**} \in \cM\right\}\\
    &\subset T_{p}(\cM),
\end{align}
In fact 
\begin{align}
    &\frac{p^{**}(X) - p^{*}(X)}{p(X)}\\
    &=\left(\frac{p^{**}(X)}{p(X)} - 1\right) - \left(\frac{p^{*}(X)}{p(X)} - 1\right) \in T_{p}(\cM).
\end{align}
Hence the conclusion.
\subsection*{Proof of Corollaries \ref{cor:cor1} and \ref{cor:cor2}}
For any $p \in \cM$, any $p^* \in \cM_{\theta, \gamma_j}(p)$ $(j = 1, 2)$ and any $s(X) \in T_{p^*}(\cM_{\theta, \gamma_j}(p))$, we have
\begin{align}
    &\langle\phi_p(X), s(X)\rangle_{p^*}\\
    &=\left\langle\prod^e_{p^* \to p}\phi_p(X), \prod^m_{p^* \to p} s(X)\right\rangle_{p}\\
    &=\left\langle\phi_p(X), \prod^m_{p^* \to p} s(X)\right\rangle_{p}\\
    &~~-\left\langle\E\{\phi_p(X)\}, \prod^m_{p^* \to p} s(X)\right\rangle_{p}\\
    &=\left\langle\phi_p(X), \prod^m_{p^* \to p} s(X)\right\rangle_{p}=0.
\end{align}
Hence, if $\phi(X; \theta, \gamma_1, \gamma_2)$ is doubly robust, then by Theorem \ref{thm:question1}, $\left\langle\phi_p(X), \prod^m_{p^* \to p} s(X)\right\rangle_{p} = 0$. This leads to Corollary \ref{cor:cor1}. On the other hand, if $\left\langle\phi_p(X), \prod^m_{p^* \to p} s(X)\right\rangle_{p} = 0$, then by Theorem \ref{thm:iff}, $\phi(X; \theta, \gamma_1, \gamma_2)$ is doubly robust. This leads to Corollary \ref{cor:cor2}.



\bibliographystyle{imsart-nameyear} 
\bibliography{ref}       


\end{document}